\documentclass[12pt]{article}

\usepackage[cp1251]{inputenc}  
\usepackage[T2A]{fontenc}      
\usepackage[russian]{babel}    

\makeatletter
\@addtoreset{equation}{section}
\makeatother

\begin{document}
{
 \sloppy

\begin{center}
{\bf \Large\bf On the measuring of independence degree of the two discrete random variables.}
\\ \bigskip
E.А.Yanovich
\\ \bigskip
University of St. Petersburg,\\
  Department of Higher Mathematics,\\
  Moskovskii 151a, 59, St. Petersburg, 196128 Russia\\
   E-mail: teduard@land.ru
\end{center}

\begin{quote}
{\footnotesize In this paper we construct the new coefficient which allows to measure quantitatively the independence of the two discrete random variables. The new inequalities for the matrices with non-negative elements are found.}
\end{quote}
\bigskip

\section{Introduction.}

Statistical description of the two random variables with finite number of states is determined by $(n\times m)$ matrix $P$
($n\ge2\,,m\ge2$). Its elements $p_{i,j}$ ($p_{i,j}\ge 0$) define the probabilities of simultaneous appearance $i$-th value of one random variable and $j$-th value of another one. The condition of probability preservation
\begin{equation}
\label{total_prob}
\sum_{i,j}p_{i,j}=1
\end{equation}
is valid.

Random variables will be independent if and only if the matrix $P$ elements are presented in the form $p_{i,j}=p_i\,q_j$. Then the matrix $P$ can be represented as a product of the matrix-column $B=(p_1,\ldots,p_n)^T$ and the matrix-line $C=(q_1,\ldots,q_m)$ and the condition of probability preservation~(\ref{total_prob}) take the form
$$
\sum_{i=1}^n p_i=1\,,\quad \sum_{i=1}^m q_i=1
$$

If the matrix $P$ is known then the matrices $B$ and $C$ are defined by
$$
B=P(D_m)^T\,\quad C=D_nP\,,
$$
where $D_k=(1,\ldots,1)$ is the unit line of length $k$.

Random variables is said to be completely dependent if and only if to every value of each random variable corresponds one and only one value of another random variable. That is each line and each column of $P$ contains at most one non-zero element. In this case the dependence between random variables is said to be functional.

For the measuring independence degree it is necessary to construct such number coefficient $k$ that satisfies to some natural conditions. We claim that~\cite{1}

$1.\quad k\in [0,1]$

$2.\quad\left(k=0\right)\Leftrightarrow\left(independence\right)$

$3.\quad\left(k=1\right)\Leftrightarrow\left(complete\,dependence\right)$\\
The coefficient $k$ satisfying to all three conditions can be constructed as follows
$$
k=\frac{\mu}{\mu_f}\,,
$$
where
$$
\mu=\sum_k(M_k)^2\,,\quad \mu_f=\sum_{1\le i<j\le n}(s_i)^2(s_j)^2\,,\quad s_i=\sum_{k=1}^mp_{i,k}
$$

Here $\mu$ is the sum of squares of all second order determinants $M_k$ of the matrix $P$. We suppose that $s_i>0$ for all $i$ and $n\le m$ (in the case $n>m$, we should transpose the matrix $P$).

\section{Proof of the 1-3 conditions.}

Let's start from the condition 2. The condition $p_{i,j}=p_i\,q_j$ is necessary and sufficient for $rank\,P=1$~\cite{2}. It follows that $k=0$ ($\mu=0$) if and only if the random variables are independent.

It is evident that $\mu\ge 0$. For the proof of the condition 1 and 3, let's consider the matrix with two lines
$$
P=\left(\begin{array}{cccc}
                              a_1&a_2&\ldots&a_m\\
                              b_1&b_2&\ldots&b_m
           \end{array}\right)\,,\quad m\ge 2
$$
We have
$$
\mu=\sum_k(M_k)^2=\sum_{i<j}|a_ib_j-a_jb_i|^2=\sum_ia_i^2\sum_ib_i^2-\left(\sum_ia_ib_i\right)^2\le
\sum_ia_i^2\sum_ib_i^2\le
$$
$$
\le(a_1+\ldots+a_m)^2\,(b_1+\ldots+b_m)^2=(s_1)^2\,(s_2)^2=\mu_f
$$

The last inequality will be exact equality if and only if each sequence $a_n$ and $b_n$ contains at most one non-zero element. On the other hand, if $s_1=a_i$, $s_2=b_j$ and $i\ne j$ then $\mu=\mu_f$. The conditions 1 and 3 are proved for $n=2$.

Let's consider a general case. Suppose that $n\le m$. Estimating the sums corresponding to each pare of lines as before, we obtain
$$
\mu\le (s_1)^2((s_2)^2+\ldots+(s_n)^2)+(s_2)^2((s_3)^2+\ldots+(s_n)^2)+\ldots+(s_{n-1})^2(s_n)^2=\mu_f\,,
$$
and $\mu=\mu_f$ if and only if each line contains one non-zero element and these elements are situated in the different columns of $P$. It is possible since $n\le m$.

Thus all conditions of the coefficient $k$ are proved.

One can note that $\mu_f$ is the maximum value of $\mu$ for fixed $s_1,s_2,\ldots,s_n$. This value should be achieved if $s_i=p_{i,k_i}$ and all $k_i$ were different that is if the dependence between random variables was functional.

Let's note that due to the condition~(\ref{total_prob}) all results remain valid and for infinite probability matrices $P$.


\begin{thebibliography}{2}
\bibitem{1}
Yanovich-Tur E.A., $"$On independence degree of the two random variables$"$, Thesis of the conference $"$Fundamental investigations and innovations in technical universities$"$, St-Petersburg, 2007.
\bibitem{2}
Faddeev D.K., Sominskii I.S., $"$Collection of tasks on higher algebra$"$, Мoskow, 1972.
\end{thebibliography}
\end{document}